\documentclass{amsart}

\usepackage[leqno]{amsmath}
\usepackage{latexsym,amsthm,amsfonts,amssymb}
\usepackage{hyperref}

\newcommand{\numberseries}{\mdseries}   
\newlength{\thmtopspace}                
\newlength{\thmbotspace}                
\newlength{\thmheadspace}               
\newlength{\thmindent}                  
\newtheoremstyle{fixed bf head,slanted body}
                {\thmtopspace}{\thmbotspace}{\slshape}
                {\thmindent}{\bfseries}{.}{\thmheadspace}
                {{\numberseries \thmnumber{(#2) }}\thmname{#1}\thmnote{ (#3)}}
\newtheoremstyle{fixed bf head,upright body}
                {\thmtopspace}{\thmbotspace}{\upshape}
                {\thmindent}{\bfseries}{.}{\thmheadspace}
                {{\numberseries \thmnumber{(#2) }}\thmname{#1}\thmnote{ (#3)}}
\newtheoremstyle{numbered paragraph}
                {\thmtopspace}{\thmbotspace}{\upshape}
                {\thmindent}{\upshape}{}{0pt}
                {{\numberseries \thmnumber{(#2) }}}
\theoremstyle{fixed bf head,slanted body}
\newtheorem{thm}{Theorem}[section]          \newtheorem*{thm*}{Theorem}
\newtheorem{prp}[thm]{Proposition}      \newtheorem*{prp*}{Proposition}
\newtheorem{lem}[thm]{Lemma}            \newtheorem*{lem*}{Lemma}
\theoremstyle{fixed bf head,upright body}
\newtheorem{exa}[thm]{Example}
\newtheorem{rmk}[thm]{Remark}
\theoremstyle{numbered paragraph}
\newtheorem{ipg}[thm]{}
\newlength{\thmlistleft}        
\newlength{\thmlistright}       
\newlength{\thmlistpartopsep}   
\newlength{\thmlisttopsep}      
\newlength{\thmlistparsep}      
\newlength{\thmlistitemsep}     
\newcounter{prt}
\newenvironment{prt}{\begin{list}{\upshape (\alph{prt})}%
    {\usecounter{prt}%
      \setlength{\leftmargin}{\thmlistleft}%
      \setlength{\labelwidth}{\thmlistleft}%
      \setlength{\rightmargin}{\thmlistright}%
      \setlength{\partopsep}{\thmlistpartopsep}%
      \setlength{\topsep}{\thmlisttopsep}%
      \setlength{\parsep}{\thmlistparsep}%
      \setlength{\itemsep}{\thmlistitemsep}}}%
  {\end{list}}%
\newcounter{eqc} 
\newenvironment{eqc}{\begin{list}{\upshape (\textit{\roman{eqc}})}%
    {\usecounter{eqc}%
      \setlength{\leftmargin}{\thmlistleft}%
      \setlength{\labelwidth}{\thmlistleft}%
      \setlength{\rightmargin}{\thmlistright}%
      \setlength{\partopsep}{\thmlistpartopsep}%
      \setlength{\topsep}{\thmlisttopsep}%
      \setlength{\parsep}{\thmlistparsep}%
      \setlength{\itemsep}{\thmlistitemsep}}}%
  {\end{list}}%
\newcommand{\eqclbl}[1]{{\upshape(\textit{#1})}}
\newenvironment{prf*}[1][Proof]{%
  \begin{proof}[\bf #1]
    \setcounter{equation}{0}
    }
  {\end{proof}
}
\newcommand{\proofofimp}[3][:]{\mbox{\eqclbl{#2}$\!\implies\!$\eqclbl{#3}#1}}
\newcommand{\pgref}[1]{(\ref{#1})}
\newcommand{\thmref}[2][Theorem~]{#1\pgref{thm:#2}}
\newcommand{\prpref}[2][Proposition~]{#1\pgref{prp:#2}}
\newcommand{\lemref}[2][Lemma~]{#1\pgref{lem:#2}}
\newcommand{\exaref}[2][Example~]{#1\pgref{exa:#2}}
\newcommand{\rmkref}[2][Remark~]{#1\pgref{rmk:#2}}
\newcommand{\secref}[2][Section~]{#1\ref{sec:#2}}
\renewcommand{\eqref}[1]{\pgref{eq:#1}}
\newcommand{\thmcite}[2][?]{\cite[thm.~#1]{#2}}

\newcommand{\prpcite}[2][?]{\cite[prop.~#1]{#2}}
\newcommand{\lemcite}[2][?]{\cite[lem.~#1]{#2}}
\newcommand{\seccite}[2][?]{\cite[sec.~#1]{#2}}
\newcommand{\rmkcite}[2][?]{\cite[rmk.~#1]{#2}}

\def\urltilda{\kern -.15em\lower .7ex\hbox{\~{}}\kern .04em} 
\newcommand{\Poin}[2][R]{\operatorname{P}^{#1}_{#2}(z)}
\newcommand{\Bass}[2][R]{\operatorname{I}_{#1}^{#2}(t)}
\newcommand{\Kz}{\operatorname{K}} 
\newcommand{\Hilb}[1]{\operatorname{H}_{#1}(z)}
\newcommand{\e}[1]{\varepsilon_{#1}}
\newcommand{\set}[2][\mspace{1mu}]{\{#1 #2 #1\}}
\newcommand{\xra}[2][]{\xrightarrow[#1]{\;#2\;}}
\newcommand{\Hom}[3][R]{\operatorname{Hom}_{#1}(#2,#3)}
\newcommand{\m}{\mathfrak{m}}
\newcommand{\n}{\mathfrak{n}}
\newcommand{\ft}{\mathfrak{t}}
\newcommand{\pows}[2][k]{#1[\mspace{-2.3mu}[#2]\mspace{-2.3mu}]}
\newcommand{\type}[2][R]{\operatorname{type}_{#1}#2}
\newcommand{\qtext}[1]{\quad\text{#1}\quad}
\newcommand{\qqtext}[1]{\qquad\text{#1}\qquad}

\newcommand{\qqand}{\qqtext{and}}
\renewcommand{\dim}[2][R]{\operatorname{dim}_{#1}#2}
\newcommand{\Tor}[4][R]{\operatorname{Tor}^{#1}_{#2}(#3,#4)}
\newcommand{\Ext}[4][R]{\operatorname{Ext}_{#1}^{#2}(#3,#4)}
\newcommand{\is}{\cong}
\renewcommand{\H}[2][]{\operatorname{H}_{#1}(#2)}
\newcommand{\deq}{\:=\:}

\newcommand{\dle}{\:\le\:}

\renewcommand{\a}{\alpha}
\renewcommand{\b}{\beta}
\renewcommand{\le}{\leqslant}
\renewcommand{\ge}{\geqslant}

\numberwithin{equation}{thm}

\begin{document}

\title[The Golod property of powers of the maximal ideal]{The Golod
  property of powers\\of the maximal ideal of a local ring}

\author[L.\,W. Christensen]{Lars Winther Christensen}

\address{Texas Tech University, Lubbock, TX 79409, U.S.A.}

\email{lars.w.christensen@ttu.edu}

\urladdr{http://www.math.ttu.edu/\urltilda lchriste}

\author[O. Veliche]{Oana Veliche}

\address{Northeastern University, Boston, MA~02115, U.S.A.}

\email{o.veliche@northeastern.edu}

\thanks{L.W.C.\ was partly supported by NSA grant H98230-14-0140 and
  Simons Foundation collaboration grant 428308.}

\date{16 December 2017} 



\begin{abstract}
  We identify minimal cases in which a power $\m^i\not=0$ of the
  maximal ideal of a local ring $R$ is not Golod, i.e.\ the
  quotient ring $R/\m^i$ is not Golod. Complementary to a 2014
  result by Rossi and \c{S}ega, we prove that for a generic artinian
  Gorenstein local ring with $\m^4=0\ne \m^3$, the quotient $R/\m^3$
  is not Golod.  This is provided that $\m$ is minimally generated by
  at least $3$ elements. Indeed, we show that if $\m$ is
  $2$-generated, then every power $\m^i\ne 0$ is Golod.
\end{abstract}

\maketitle

\thispagestyle{empty}

\section{Introduction}
\label{sec:Introduction}

\noindent
In this paper a local ring is a commutative noetherian ring $R$ with
unique maximal ideal~$\m$. Such a ring is called \emph{Golod} if the
ranks of the modules in the minimal free resolution of the residue
field $R/\m$ attain the upper bound established by Serre; the precise
definition is recalled in \secref{gor}.

The field $R/\m$ is trivially Golod, and so is~the quotient
ring $R/\m^2$; see for example Avramov's exposition
\prpcite[5.2.4]{ifr}. Moreover, if $R$ is a regular local ring, then
the quotient $R/\m^i$ is Golod for every $i\ge 1$. Rossi and \c{S}ega
\cite{MERLMS14} prove that for a generic artinian Gorenstein local
ring $(R,\m)$ with $\m^4 \ne 0$, every proper quotient $R/\m^i$ is
Golod.

In this note we provide minimal examples of local rings with proper
quotients $R/\m^i$ that are \emph{not} Golod. They come out of an
investigation of the complementary case to above mentioned result from
\cite{MERLMS14}.  The following extract from \thmref{T} points to a
whole family of local rings with $\m^4=0$ and $R/\m^3$ not Golod. In
fact, this is the behavior of generic graded Gorenstein local
$k$-algebras of socle degree $3$.

\begin{thm}
  \label{thm:t}
  Let $k$ be a field; set $Q = \pows{x,y,z}$ and $\n = (x,y,z)$. Let
  $I\subseteq \n^2$ be a homogeneous Gorenstein ideal in $Q$ with $\n^4 \subseteq I
  \not\supseteq \n^3$ and set $(R,\m) = (Q/I, \n/I)$. The following
  conditions are equivalent.
  \begin{eqc}
  \item $I$ is generated by quadratic forms.
  \item $R$ is Koszul, i.e.\ the minimal free resolution of $k$ over
    $R$ is linear.
  \item $R$ is complete intersection.
  \item $R$ has an exact zero divisor, i.e.\ an element $a \ne 0$ with
    $(0:a)$ principal.
  \item $R/\m^3$ is not Golod.
  \end{eqc}
\end{thm}
\noindent
To discuss in which sense these rings are generic and constitute
minimal examples of local rings $(R,\m)$ with proper quotients
$R/\m^i$ that are not Golod, we start to introduce the terminology
that will be used throughout the paper.
\begin{equation*}
  \ast \ \ast \ \ast
\end{equation*}
Let $(R,\m)$ be a local ring with the residue field $k=R/\m$.

\begin{ipg}
  \label{i} As remarked above, the quotient rings $R/\m$ and $R/\m^2$
  are both Golod.

  Thus, for a proper quotient $R/\m^i$ to be not Golod, one must have
  $i\ge 3$.
\end{ipg}

\begin{ipg}
  \label{e}
  For an ideal $I \subseteq R$ we denote by $\mu(I)$ its minimal
  number of generators. The number $\mu(\m)$ is the \emph{embedding
    dimension} of $R$.  It is known that every local ring of embedding
  dimension $1$ is Golod; see \pgref{c}. Further, we prove in
  \thmref{e2} that if $R$ has embedding dimension $2$, then every
  proper quotient $R/\m^i$ is Golod.

  Thus, for a proper quotient $R/\m^i$ to be not Golod, the
  embedding dimension of $R$ must be at least 3.
\end{ipg}

\begin{ipg}
  \label{c} The embedding dimension $e$ and the Krull dimension $d$ of
  $R$ satisfy $e\ge d$, with equality if and only if $R$ is
  regular. The difference $e-d$ is called the \emph{codimension} of
  $R$.  If $R$ has codimension at most $1$, then $R$ is Golod, see
  \prpcite[5.2.5]{ifr}, and $R/\m^i$ is known to be Golod for every
  $i\ge 1$ by work of \c{S}ega \prpcite[6.10]{LMS01}.

  Thus, for a proper quotient $R/\m^i$ to be not Golod, the
  codimension of $R$ must be at least $2$.
\end{ipg}

\noindent
We exhibit in \exaref{e2} a complete intersection local ring $(R,\m)$
of embedding dimension 3 and codimension 2 with $R/\m^i$ not Golod for
all $i\ge 3$. Thus, among local rings $(R,\m)$ with the property that
a proper quotient $R/\m^i$ is not Golod, this ring is minimal with
regard to codimension. It has dimension $1$; a $0$-dimensional, i.e.\
artinian, local ring with the property must have codimension at least
$3$.

\begin{ipg}
  \label{a}
  Let $(R,\m)$ be artinian. The integer $s$ with $\m^{s+1}=0\not=\m^s$
  is called the \emph{socle degree} of $R$. The \emph{type}
  of $R$ is the dimension of the socle as a vector space over the residue
  field; i.e.\ $\type[]{R} = \dim[k]{(0:\m)}$.

  In view of \pgref{i}, the ring $R$ must have socle degree at least
  $3$ for a proper quotient $R/\m^i$ to be not Golod.
\end{ipg}

It follows that among local rings $(R,\m)$ with a proper quotient
$R/\m^i$ that is not Golod, the rings in \thmref{t} are minimal with
regard to dimension and embedding dimension, and then with regard to
socle degree and type.

\begin{ipg}
  An element $a\ne 0$ in a commutative ring $R$ is called an
  \emph{exact zero divisor} if the annihilator $(0 : a)$ is a
  principal ideal. Exact zero divisors in artinian Gorenstein local
  rings of socle degree $3$ were studied in depth by Henriques and
  \c{S}ega in \cite{IBHLMS11}.
\end{ipg}

A generic artinian Gorenstein local graded $k$-algebra of socle degree
$3$ has an exact zero divisor; this follows from work of Conca, Rossi,
and Valla \cite{CRV-01}; see \rmkcite[4.3]{IBHLMS11}\footnote{The
  notion of \emph{generic} that is used in \cite{CRV-01} and \cite{IBHLMS11}
  is, at least formally, different from the one used in
  \cite{MERLMS14}. However, \thmref{T} and thereby \thmref{t} apply
  to rings that are generic in either sense; see also the discussion
  after \pgref{MCP}.}. In particular, generic Gorenstein algebras of
embedding dimension $3$ and socle degree $3$ satisfy the conditions in
\thmref{t}.  We provide examples of Gorenstein algebras of socle
degree $3$ without exact zero divisors in \prpref[]{E} and
\exaref[]{C}.

\section{Embedding dimension $2$}
\label{sec:e2}

\noindent
We prove that all proper quotients $R/\m^i$ are Golod for any local
ring $(R,\m)$ of embedding dimension $2$. To frame this result we
exhibit a local ring of embedding dimension $3$ and codimension $2$
such that $R/\m^i$ is not Golod for $i\ge 3$.

\begin{ipg}
  Let $(R,\m)$ be a local ring. The \emph{valuation} of an ideal $J
  \subset R$ is the largest integer $i$ with $J \subseteq \m^i$; it is
  written $v_R(J)$. For an element $x\in R$ one sets $v_R(x) =
  v_R((x))$.
\end{ipg}

In \cite{GSc64} Scheja shows that a local ring of embedding dimension
$2$ is either complete intersection or Golod; see also
\prpcite[(5.3.4)]{ifr}. The gist of the next result is that such a
complete intersection cannot arise as a proper quotient by a power
of the maximal ideal.

\begin{thm}
  \label{thm:e2}
  Let $(R,\m)$ be a local ring of embedding dimension $2$. For every
  $i \ge 1$ with $\m^i\ne 0$, the quotient ring $R/\m^i$ is Golod.
\end{thm}

\begin{prf*}
  Let $(\widehat{R},\widehat{\m})$ be the $\m$-adic completion of
  $R$. One has $R/\m^i \is \widehat{R}/\widehat{\m}^i$, so we may
  assume that $R$ is complete. By Cohen's Structure Theorem there is a
  regular local ring $(Q,\n)$ of embedding dimension $2$ and an ideal
  $J \subseteq \n^2$ with $R \is Q/J$. Thus, one has $R/\m^i \is
  Q/(J+\n^i)$, and $\m^i \ne 0$ if and only if $\n^i \not\subseteq
  J$. Per \pgref{i} we may assume that one has $i\ge
  3$. To prove that $R/\m^i$ is Golod, it suffices to argue that
  it is not complete intersection. Thus,
  we now argue that if $\n^i \not\subseteq J$, then $J + \n^i$ cannot
  be generated by two elements.

  For every $i$ the ideal $\n^i$ is minimally generated by $i+1$
  elements. In particular, if $J \subseteq \n^i$, then $J+\n^i$ cannot
  be generated by $2$ elements. We now assume that $J$ is not
  contained in $\n^i$; setting $t = v_Q(J)$ this means $t+1 \le
  i$. Write
  \begin{equation*}
    J = (f_1,\ldots,f_n)
  \end{equation*}
  where $v_Q(f_1)=t$, and assume towards a contradiction that $J +
  \n^i$ is generated by $2$ elements. As $\n(J+\n^i) \subseteq
  \n^{t+1}$ and $f_1 \in \n^t\setminus \n^{t+1}$ we may assume that $J
  + \n^i$ is minimally generated by $f_1$ and some element $g$. One
  has $f_n = af_1 + bg$ for $a,b\in Q$, and $b$ is not a unit as
  $g\not\in J$. Now write
  \begin{equation*}
    g = \sum_{j=1}^na_jf_j + h
  \end{equation*}
  with $h \in \n^i$; without loss of generality we may assume $a_1
  =0$.  Now one has
  \begin{equation*}
    f_n = af_1 + b\bigg(\sum_{j=2}^na_jf_j + h\bigg)\:,
  \end{equation*}
  whence $f_n(1-a_nb)$ belongs to $(f_1,\ldots,f_{n-1}) + \n^i$. As
  $1-a_nb$ is a unit, this yields
  \begin{equation*}
    J + \n^i = (f_1,\ldots,f_{n-1}) +  \n^i\;.
  \end{equation*}
  By recursion, one gets $J + \n^i = (f_1) + \n^i$, and it follows
  from \lemref{e2} that $J + \n^i$ cannot be generated by 2 elements;
  a contradiction.
\end{prf*}

\begin{lem}
  \label{lem:e2}
  Let $(Q,\n)$ be a regular local ring of embedding dimension $2$. If
  $f\in\n^2$ and $i\ge 2$ are such that $\n^i \not\subseteq (f)$, then
  one has $\mu((f) + \n^i) > 2$.
\end{lem}

\begin{prf*}
  If $ f\in \n^i$, then the statement is clear as one has $\mu(\n^i) =
  i+1 \ge 3$; thus we may assume that $f \not\in\n^i$.  Assume towards
  a contradiction that $(f) + \n^i$ is minimally generated by two
  elements. Set $t = v_Q(f)$; as one has $\n((f)+\n^i) \subseteq
  \n^{t+1}$ and $f \in \n^t\setminus \n^{t+1}$ we may assume that $(f)
  + \n^i$ is minimally generated by $f$ and some element $g$ of
  valuation $i$.  Let $x$ and $y$ be minimal generators of $\n$ and
  write
  \begin{equation*}
    g = \sum_{j=0}^ia_jx^{i-j}y^j \qqand x^{i-j}y^j = b_jf + c_jg
  \end{equation*}
  with $a_j,b_j,c_j \in Q$. These expressions yield $g =
  \big(\sum_{j=0}^ia_jb_j\big)f + \big(\sum_{j=0}^ia_jc_j\big)g$, so
  if $c_j\in\n$ for all $j$ one has $g =
  \big(1-\sum_{j=0}^ia_jc_j\big)^{-1}\big(\sum_{j=0}^ia_jb_j\big)f$,
  which contradicts the assumption that $f$ and $g$ minimally generate
  $(f)+\n^i$.  Thus, $c_j$ is a unit for some $j$, whence $(f) + \n^i$ is
  minimally generated by $f$ and $x^{i-j}y^j$. By symmetry in $x$ and
  $y$ we may assume that $j \ge 1$. Write $x^{i-j+1}y^{j-1} = af +
  bx^{i-j}y^j$ with $a,b\in Q$; one then has
  \begin{equation*}
    af = x^{i-j}y^{j-1}(x-by)\:.
  \end{equation*}

  First consider the case $j=i$ and recall that $i\ge 2$. Notice that
  $f \in (y)$ would force $(f) + \n^i = (f,y^i) \subseteq (y)$, which
  is absurd, so $f$ is not divisible by $y$. Since $Q$ is a domain and
  $(y)$ is a prime ideal, it follows that $a$ is divisible by
  $y^{j-1}$. Writing $a = a'y^{j-1}$ one now has $a'f = x - by$, which
  is absurd as $f\in\n^2$ and $x,y$ minimally generate $\n$.

  Finally consider the case $j < i$. Notice as above that $f$ is not
  divisible by $x$ nor by $y$. Therefore, one has $a =
  a'x^{i-j}y^{j-1}$ and $a'f = x-by$, which again contradicts the
  assumption that $x$ and $y$ minimally generate $\n$.
\end{prf*}

Apart from showing that quotients $R/\m^i$ of a codimension $2$
local ring may not be Golod, the next example shows that there are
artinian local rings of embedding dimension $3$ and any socle degree
$s\ge 3$ with $R/\m^i$ not Golod for $3 \le i \le s$.

\begin{exa}
  \label{exa:e2}
  Let $k$ be a field and consider the codimension 2 complete
  intersection local ring $R = \pows{x,y,z}/(x^2,y^2)$ with maximal
  ideal $\m = (x,y,z)/(x^2,y^2)$. For $i \ge 3$ it is elementary to
  verify that the quotient ring
  \begin{equation*}
    R/\m^i = \frac{\pows{x,y,z}}{(x^2,y^2,z^i,xz^{i-1},yz^{i-1},xyz^{i-2})}
  \end{equation*}
  is not Golod. Indeed, the Koszul complex $K$ on the generators
  ${x,y,z}$ of $\m$ is the exterior algebra of the free module with
  basis $\e{x}, \e{y}, \e{z}$. The element $xy(\e{x}\wedge\e{y})$ in
  $K_2$ is the product of the cycles $x\e{x},y\e{y}\in K_1$, and it is
  not a boundary as one has
  \begin{equation*}
    \partial(w(\e{x}\wedge\e{y}\wedge\e{z})) = xw(\e{y}\wedge\e{z}) -
    yw(\e{x}\wedge\e{z}) + zw(\e{x}\wedge\e{y})\:,
  \end{equation*}
  and $xy \not\in (z)$. Thus there is a non-trivial product in Koszul
  homology, which by Golod's original work \cite{ESG62} means that $R/
  \m^i$ is not Golod.

  Notice that for $s \ge 3$ the artinian local ring $(\dot{R},\dot{\m}) =
  (R/\m^{s+1},\m/\m^{s+1})$ has $\dot{R}/\dot{\m}^i = R/\m^i$ not
  Golod for $3 \le i \le s$.
\end{exa}

\section{Artinian Gorenstein rings}
\label{sec:gor}

\noindent
Let $(R,\m,k)$ be a local ring, i.e.\ $k$ is the residue field
$R/\m$. For a finitely generated $R$-module $M$, the power series
$\Poin{M} = \sum_{j=0}^\infty (\dim[k]{\Tor{j}{k}{M}})z^j$ is called
the \emph{Poincar\'e series} of $M$ over $R$.  The numbers
$\dim[k]{\Tor{*}{k}{M}}$ are known as the \emph{Betti numbers} of $M$;
they record the ranks of the free modules in the minimal free
resolution of $M$ over $R$. In particular, $R$ is Golod if and only if
one has
\begin{equation}
  \label{eq:golod}
  \Poin{k} = \frac{(1+z)^e}
  {1 - \sum_{j=1}^{e-d} (\dim[k]{\H[j]{\Kz^R}})z^{j+1}}
\end{equation}
where $e$ and $d$ are the embedding dimension and depth of $R$, and
$\Kz^R$ is the Koszul complex on a minimal set of generators of $\m$;
see \cite[(5.0.1)]{ifr}.

It is standard to refer to $\Poin{k}$ as the Poincar\'e series of
$R$. For an artinian Gorenstein local ring of embedding dimension
$e\ge 2$ and socle degree $s$, the Poincar\'e series of $R/\m^s$ was
computed by Avramov and Levin \thmcite[2]{GLLLLA78}:
\begin{equation}
  \label{eq:LA}
  \Poin[R/\m^s]{k} = \frac{\Poin{k}}{1-z^2\Poin{k}}\;.
\end{equation}
The special case of a complete intersection was first done by
Gulliksen \thmcite[1]{THG72}:
\begin{equation}
  \label{eq:THG}
  \Poin[R/\m^s]{k} = \frac{1}{(1-z)^e - z^2}\;.
\end{equation}

\begin{prp}
  \label{prp:ci}
  If $(R,\m)$ is an artinian complete intersection local ring of
  embedding dimension at least $3$ and socle degree $s$, then $R/\m^s$ is
  not a Golod ring.
\end{prp}

\begin{prf*}
  Let $e$ be the embedding dimension of $R$. By \eqref{THG} one has
  \begin{equation*}
    \Poin[R/\m^s]{k} = \frac{1}{(1-z)^e-z^2}
    = \frac{(1+z)^e}{(1-z^2)^e - z^2(1+z)^e}\:;
  \end{equation*}
  notice that the denominator has degree $2e$ since $e\ge 3$.
  If $R/\m^s$ were Golod, then by \eqref{golod} the
  Poincar\'e series would have the form $(1+z)^e/d(z)$ where the
  denominator $d(z)$ has degree $e+1 < 2e$. Thus,  $R/\m^s$ is
  not a Golod ring.
\end{prf*}

\begin{ipg}
  \label{MCP}
  Let $(R,\m,k)$ be artinian of embedding dimension $e$ and socle
  degree $s$.

  By Cohen's Structure Theorem there is a regular local ring $(Q,\n)$
  and an ideal $I$ with $\n^{s+1} \subseteq I \subseteq \n^2$ such
  that $Q/I \is R$. This is called the \emph{minimal Cohen
    presentation} of $R$; notice in particular that $Q$ also has
  embedding dimension $e$.

  Denote by $h_R(i)$ and $\Hilb{R}$ the Hilbert function and Hilbert
  series of $R$; i.e.\
  \begin{equation*}
    h_R(i) = \dim[k](\m^i/\m^{i+1}) \ \text{for } i\ge 0 
    \qqand \Hilb{R} = \sum_{i=0}^s h_R(i)z^i\:.
  \end{equation*}
  One says that $R$ is \emph{Koszul} if the associated graded
  $k$-algebra $\bigoplus_{i=0}^{s} \m^i/\m^{i+1}$ is Koszul in the
  traditional sense that $k$ has a linear resolution; see the
  discussion in \cite[1.10]{IBHLMS11}.  If $R$ is Koszul, then one has
  $\Poin{k}\H[R]{-z} = 1$.

  If $R$ is Gorenstein, then for every $i\ge 0$ there is an
  inequality
  \begin{equation*}
    h_R(i) \dle \min\set{h_Q(i),h_Q(s-i)} \deq 
    \min\left\{ {e-1+i\choose e-1}, {e-1+s-i\choose e-1} \right\}\:;
  \end{equation*}
  If equality holds for every $i$, then $R$ is called
  \emph{compressed}, see \seccite[4]{MERLMS14}.
\end{ipg}

The idea of compressed rings was introduced by Iarrobino, in
\thmcite[I]{AIr84} he shows that generic artinian Gorenstein local
standard graded algebras over a field are compressed. Rossi and
\c{S}ega prove \prpcite[6.3]{MERLMS14} that for a compressed artinian
Gorenstein local ring $(R,\m)$ of socle degree $s \ne 3$ the quotient
ring $R/\m^i$ is Golod for all $1 \le i \le s$.

Here we focus on rings of socle degree $3$. Our main result,
\thmref{T}, is a simultaneous converse to \prpref[Propositions~]{ci}
and \prpref[]{ezd} in embedding dimension $3$.

\begin{prp}
  \label{prp:ezd}
  Let $(R,\m)$ be an artinian Gorenstein local ring of embedding
  dimension at least $3$ and socle degree $3$. If $R$ has an exact
  zero divisor, then $R$ is compressed and Koszul, and the quotient
  ring $R/\m^3$ is not Golod.
\end{prp}

\begin{prf*}
  Let $e$ denote the embedding dimension of $R$.  By
  \thmcite[3.3 and prop.~4.1]{IBHLMS11} the existence of an exact zero divisor
  implies that $R$ is Koszul with Hilbert series $1 + ez + ez^2 + z^3$; hence
  $R$ is compressed.  Further one has
  \begin{equation*}
    \Poin{k} =
    \frac{1}{1-ez +ez^2-z^3} =
    \frac{(1+z)^e}{(1+z)^e(1-ez +ez^2-z^3)}\:.
  \end{equation*}
  Note that the denominator $(1+z)^e(1-ez +ez^2-z^3)$ is a polynomial
  of degree $e+3$.  Let $Q$ be the regular ring of a minimal Cohen
  presentation of $R$; cf.~\pgref{MCP}. By \prpcite[6.2]{MERLMS14} the
  ring $R/\m^3$ is Golod if and only if one has
  \begin{equation*}
    \Poin{k} =
    \frac{(1+z)^e}{1-z(\Poin[Q]{R}-1)+z^{e+1}(z+1)}\:.
  \end{equation*}
  Since $\Poin[Q]{R}$ is a polynomial of degree $e$, the denominator
  above is a polynomial of degree $e+2$. Thus $R/\m^3$ is not Golod.
\end{prf*}

\noindent
To frame \prpref[Propositions~]{ci} and \prpref[]{ezd} we show how to
produce a compressed artinian Gorenstein local ring
$(T,\ft)$ of socle degree $3$ such that $T/\ft^3$ is not Golod, though
$T$ is not complete intersection and does not have an exact zero
divisor. See also \rmkref{frame}.

\begin{prp}
  \label{prp:E}
  Let $k$ be a field and $(R,\m)$ an artinian standard graded local
  $k$-algebra of embedding dimension $e \ge 3$ and socle degree
  $2$. If $R$ is not Gorenstein and admits a non-zero minimal acyclic
  complex $F$ of finitely generated free modules, then the graded
  local $k$-algebra
  \begin{equation*}
    T = R \ltimes \mathsf{\Sigma}\Hom[k]{R}{k} \qtext{with maximal ideal} 
    \ft = \m \ltimes \mathsf{\Sigma}\Hom[k]{R}{k}
  \end{equation*}
  is Gorenstein and compressed with Hilbert series $1 + (2e-1)z +
  (2e-1)z^2 + z^3$.  Furthermore, the following
  hold:
  \begin{prt}
  \item $T$ is not complete intersection.
  \item If $R$ does not have an exact zero-divisor, then neither does $T$.
  \item If\, $\H[n]{\Hom{F}{R}}=0$ holds for some $n$, then
    $T/\ft^3$ is not Golod.
  \end{prt}
\end{prp}

\noindent
For a concrete example of a local $k$-algebra $(R,\m)$ that meets the
assumptions in the Proposition---including those in parts (b) and
(c)---see Christensen, Jorgensen, Rahmati, Striuli, and Wiegand
\seccite[9]{CJRSW-12}.

\begin{prf*}
  From \thmcite[A]{LWCOVl07} it is known that the Hilbert series of
  $R$ is $1 + ez + (e-1)z^2$. As a graded $k$-vector space $T$ has the
  form $R \oplus \mathsf{\Sigma}\Hom[k]{R}{k}$, so one has
  \begin{equation*}
    \Hilb{T} = \Hilb{R} + z^3\H[R]{z^{-1}} = 1 +
    (2e-1)z + (2e-1)z^2 +z^3.
  \end{equation*}
  Recall that $E = \Hom[k]{R}{k}$ is the injective envelope of $k$
  over $R$. As a local ring, $T$ is the trivial extension of $R$ by
  $E$, so by \cite[lem.\ in sec.\ 3]{THG72} it is Gorenstein, and
  evidently it is compressed; cf.~\pgref{MCP}.  

  In the sequel, let $Q/I$ be a minimal Cohen presentation of $T$.

  (a): If $T$ were complete intersection, then $I$ would be
  generated by $2e-1$ elements, but that is not possible
  as one has
  \begin{equation*}
    h_Q(2)-h_T(2) = \binom{2e}{2} - (2e-1) = (e-1)(2e-1) > 2e-1\:.
  \end{equation*}

  (b): Assume towards a contradiction that $(x,\a)$ is an exact zero
  divisor in $T$ with annihilator generated by $(y,\b)$. The element
  $(y,\b)$ is also an exact zero divisor, called the complementary
  divisor; see \rmkcite[1.1]{IBHLMS11}. It follows from
  \prpcite[4.1]{IBHLMS11} that $(x,\a)$ and $(y,\b)$ belong to
  $\ft\setminus \ft^2$. Hence $x$ and $y$ belong to
  $\m\setminus\m^2$ and, evidently, one has $xy=0$. Any element in
  $\m^2 \ltimes 0$ annihilates $(x,\a)$ and is hence contained in the
  ideal generated by $(y,\b)$. In particular, one has $\m^2 =y\m$ and
  by symmetry $\m^2 = x\m$. Now it follows from
  \lemcite[4.3(c)]{CJRSW-12} that $x$ and $y$ are exact zero-divisors
  in $R$, a contradiction.

  (c) We argue that $T/\ft^3$ is not Golod by comparing two
  expressions for the Poincar\'{e} series of $T$.  By a computation of
  Gulliksen \thmcite[2]{THG72} one has
  \begin{equation*}
    \Poin[T]{k} =
    \frac{\Poin{k}}{1-z\Poin{E}}\:.
  \end{equation*}
  The standard isomorphisms $\Tor{*}{k}{\Hom[k]{R}{k}} \is
  \Hom[k]{\Ext{*}{k}{R}}{k}$ and \thmcite[A]{LWCOVl07} yield:
  \begin{equation*}
    \Poin{k} = \frac{1}{(1-z)(1-(e-1)z)} \qqand \Poin{E} = 
    \Bass{} = \frac{e-1-z}{1-(e-1)z}\:.
  \end{equation*}
  Finally, a direct computation yields
  \begin{equation*}
    \Poin[T]{k} = \frac{(1+z)^{2e-1}}{(1-z)(1-2(e-1)z+z^2)(1+z)^{2e-1}}\:;
  \end{equation*}
  notice that the denominator has degree $2e+2$. On the other hand,
  the regular ring $Q$ has embedding dimension $2e-1$; in particular,
  $\Poin[Q]{T}$ is a polynomial of degree $2e-1$. As in the proof of
  \prpref{ezd} above, it follows from \prpcite[6.2]{MERLMS14} that the
  ring $T/\ft^3$ is Golod if and only if $\Poin[T]{k}$ has the form
  $(1+z)^{2e-1}/d(z)$ where $d(z)$ is a polynomial of degree $2e+1$.
\end{prf*}

\section{Embedding dimension 3 and socle degree 3}

\noindent
Let $(R,\m,k)$ be a local ring of embedding dimension $3$, let $\Kz^R$
be the Koszul complex on a minimal set of generators of $\m$, and set
$A = \H[]{\Kz^R}$. The Koszul complex is a differential graded
algebra, and the product on $\Kz^R$ induces a graded-commutative
$k$-algebra structure on $A$.  As one has $A_{\ge 4} = 0$ it follows
from Golod's original work \cite{ESG62} that $R$ is Golod if and only
if $A$ has trivial multiplication, i.e.\ $A_{\ge 1}\cdot A_{\ge 1}
=0$. Moreover, it is known from work of Assmus \thmcite[2.7]{EFA59}
that $R$ is complete intersection if and only if $A$ is isomorphic to
the exterior algebra on $A_1$.

There is a complete classification, due to Weyman~\cite{JWm89} and
Avramov et al.~\cite{LLA12,AKM-88}, of artinian local rings of
embedding dimension $3$---even more generally of local rings of
codepth $\le 3$---based on multiplication in Koszul homology. For the
precise statement of our main theorem, we need to recall one more
class from this scheme: It is called $\mathbf{T}$, and if $R$ belongs
to this class one has $A_1\cdot A_1 \ne 0 = A_1\cdot A_2$; in
particular $R$ is neither Golod nor complete intersection.

\begin{rmk}
  An artinian local ring $(R,\m)$ of embedding dimension $2$ and socle
  degree $2$ is Gorenstein if and only if it is complete intersection
  if and only if it has an exact zero divisor if and only if every
  element in $\m\setminus\m^2$ is an exact zero divisor; see
  \rmkcite[(7.1)]{CJRSW-12}. Such rings have Hilbert series $1 +2z
  +z^2$, so they are~compressed. 

  In \thmref{T} and \rmkref{T} we show that much of this behavior
  extends to embedding dimension and socle degree $3$ but perhaps not
  further.
\end{rmk}

\begin{thm}
  \label{thm:T}
  Let $(R,\m)$ be an artinian Gorenstein local ring of embedding
  dimension $3$ and socle degree $3$. The following conditions are
  equivalent.
  \begin{eqc}
  \item $R$ is complete intersection.
  \item $R$ is compressed and Koszul.
  \item $R$ has an exact zero divisor.
  \item $R/\m^3$ belongs to the class $\mathbf{T}$.
  \item $R/\m^3$ is not Golod.
  \end{eqc}
\end{thm}

\noindent
To see that \thmref{t} follows from this statement, recall that a
standard graded Koszul algebra is quadratic. Further, being
Gorenstein, the ring $Q/I$ has a symmetric Hilbert series, i.e.\ it is $1 +
3z + 3z^2+z^3$. Thus, if $I$ is quadratic, then it is minimally
generated by 3 elements, which necessarily form a regular
sequence.

\begin{prf*}
  Let $Q/I$ be a minimal Cohen presentation of $R$ and set $S =
  R/\m^3$; cf.~\pgref{MCP}.

  \proofofimp{i}{iv} If $R$ is complete intersection, then one has
  \begin{equation*}
    \Poin[S]{k} = \frac{1}{1-3z+2z^2-z^3} = \frac{(1+z)^2}{1-z-3z^2-z^5}
  \end{equation*}
  by \eqref{THG}, and that identifies $S$ as belonging to the class
  \textbf{T}; see \thmcite[2.1]{LLA12}.

  \proofofimp{iv}{v} Evident as rings of class \textbf{T} are not
  Golod.

  \proofofimp[ by contraposition:]{v}{i} If $R$ is not complete
  intersection, then \thmcite[2.1]{LLA12} yields $\Poin{k} =
  (1+z)^2/g(z)$ with $g(z) = 1 - z - (\mu(I)-1)z^2 -z^3 + z^4$. By
  \eqref{LA} one then has
  \begin{equation*}
    \Poin[S]{k} = \frac{(1+z)^2}{g(z) - z^2(1+z)^2} =  
    \frac{(1+z)^2}{1 - z - \mu(I)z^2 - 3z^3}\:,
  \end{equation*}
  and, again by \emph{loc.~cit.,} that identifies $S$ as being Golod.

  \proofofimp{i}{iii} It follows from \pgref{MCP} that $R$ has length
  at most $1+3+3+1 = 8$ with equality if and only if $R$
  compressed. On the other hand, the length of $R$ is at least $2^3 =
  8$ by \cite[\S7, prop.~7]{bour89}. Thus $R$ is compressed with
  Hilbert series $1 + 3z + 3z^2 + z^3$; in particular, one has
  $(0:\m^2) = \m^2$; see \prpcite[4.2(b)]{MERLMS14}.

  An application of \lemcite[2.8]{CRV-01} to the associated graded
  ring yields an element $\ell \in \m\setminus\m^2$ with
  $\ell\m \ne \m^2$. We argue that $\ell$ is annihilated by an element
  in $\m \setminus \m^2$.  Assume towards a contradiction that
  $(0:\ell) \subseteq \m^2$ holds. As $R$ Gorenstein, one now has
  $(\ell) \supseteq (0:\m^2) =\m^2$, and therefore, $\ell\m = \m^2$ by
  \rmkcite[2.2(1)]{IBHLMS11}, which is a contradiction. Thus there
  exists an $\ell' \in \m\setminus \m^2$ with $\ell\ell'=0$. Notice
  that $\ell,\ell' \not\in \m^2 = (0:\m^2)$ implies $\ell\m^2 = \m^3 =
  \ell'\m^2$.

  To prove that $\ell$ and $\ell'$ are exact zero divisors, it now
  suffices by \prpcite[4.1]{IBHLMS11} to verify that $\mu(\ell\m) = 2
  = \mu(\ell'\m)$; equivalently, that the linear maps $\m/\m^2 \to
  \m^2/\m^3$ given by multiplication by $[\ell]_{\m^2}$ and
  $[\ell']_{\m^2}$ have kernels of rank $1$. As $\ell\ell'=0$ the
  kernels have rank at least $1$, and by symmetry it is sufficient to
  show that the rank is $1$ for multiplication by $[\ell]_{\m^2}$.

  To this end, let $[\ell'']_{\m^2}\ne 0$ be an element in the kernel
  of $\smash{l\colon \m/\m^2 \xra{[\ell]} \m^2/\m^3}$. That is, one
  has $\ell\ell''\in \m^3$ i.e.\ $\ell\ell'' = \ell u$ for some $u \in
  \m^2$.  The elements $\ell$, $\ell'$, and $\ell''$ lift to elements
  in $\n\setminus \n^2$, and $u$ lifts to an element in $\n^2$; we
  denote these lifts by the same symbols. In $Q$ one now has
  $\ell\ell' \in I$ and $\ell(\ell'' - u) \in I$. Let $f,g,h \in \n^2$
  be a regular sequence that generates $I$ and write
  \begin{equation*}
    \ell\ell' = a'f + b'g + c'h \qqand 
    \ell(\ell'' - u) = a''f + b''g + c''h
  \end{equation*}
  with coefficients $a',\ldots,c''$ in $Q$. As the associated graded
  ring of $Q$ is a domain, one has $v_Q(\ell\ell') = v_Q(\ell) +
  v_Q(\ell') = 1+1 =2$. That is, $\ell\ell'$ is in $\n^2\setminus
  \n^3$, so we may without loss of generality assume that $c'$ is a
  unit. After cross multiplication by $\ell'$ and $\ell''-u$ and
  elimination of $\ell\ell'(\ell''-u)$ one gets $((\ell''-u)c' -
  \ell'c'')h \in (f,g)$. As $f,g,h$ is a regular sequence, this
  implies $((\ell''-u)c' - \ell'c'') \in (f,g) \subseteq \n^2$, and
  since $u$ is in $\n^2$ it further implies $\ell''c' - \ell'c'' \in
  \n^2$. Recall that $c'$ is a unit. If $c''$ were in $\n$, one would
  have $\ell'' \in \n^2$, contrary to the assumptions. Thus $c''$ is a
  unit, whence $\ell'$ and $\ell''$ are linearly dependent mod
  $\n^2$. That is, $[\ell']_{\m^2}$ spans the kernel of $l$.

  \proofofimp{iii}{ii} By \prpref{ezd}.

  \proofofimp{ii}{i} One has $\Poin{k} = 1/\H[R]{-z} =
  1/(1-3z+3z^2-z^3) = 1/(1-z)^3$, which means that $R$ is complete
  intersection; see \thmcite[2.1]{LLA12}.
\end{prf*}

While \prpcite[6.3]{MERLMS14} is a statement about proper quotients
$R/\m^i$ of all compressed Gorenstein rings of socle degree not $3$,
there is no uniform behavior of those of socle degree $3$.

\begin{exa}
  \label{exa:C}
  Let $k$ be a field. By a result of Buchsbaum and Eisenbud
  \thmcite[2.1]{DABDEs77}, the defining ideal of a Gorenstein ring
  $R=k[[x,y,z]]/I$ is generated by the sub-maximal Pfaffians of an
  odd-sized skew-symmetric matrix.

  The ideal generated by the $4\times 4$ Pfaffians of the matrix
  \begin{equation*}
    \begin{pmatrix}
      0 & x+y & 0 & 0 & y\\
      -x-y & 0 & 0 & y^2+z^2 & yz\\
      0& 0& 0& x+z& z\\
      0 & -y^2-z^2 & -x-z & 0 & x\\
      -y & -yz & -z& -x & 0
    \end{pmatrix}
  \end{equation*}
  has the form
  \begin{gather*}
    I = (xz + yz, xy +yz,x^2-yz,yz^2+z^3,y^3-z^3)\:.
  \end{gather*}
  It is straightforward to verify that a graded basis for $R$ is
  \begin{equation*}
    1\:;\quad x,y,z\:;\quad y^2,yz,z^2\:;\quad z^3\;.
  \end{equation*}
  Thus $R$ is a compressed artinian Gorenstein ring of socle degree
  3. It is not complete intersection, so by \thmref{T} it does not
  have exact zero-divisors.
\end{exa}

\section{Remarks on embedding dimension $4$}
\label{sec:4}

\noindent
The proof of the implication $(i) \Rightarrow (iii)$ in \thmref{T} relies
on $R$ being compressed, so it seems fitting to record the following remark.

\begin{rmk}
  \label{rmk:T}
  A compressed artinian
  Gorenstein ring $(R,\m)$ of embedding dimension $e \ge 4$ cannot be
  complete intersection. Indeed, let $Q/I$ be a minimal Cohen
  presentation of $R$ and let $s$ denote the socle degree. As
   $R$ is compressed, the initial degree of $I$ is $t =
  \min_i\set{h_Q(s-i) < h_Q(i)}$, and by \prpcite[4.2]{MERLMS14} one has
  $t=\lceil\frac{s+1}{2}\rceil$. A straightforward computation
  yields
  \begin{equation*}
    \mu(I) \ge h_Q(t) - h_Q(s-t) = 
    \begin{cases}
      \binom{e-2+t}{e-2} &\text{ for odd $s$}\\
      \binom{e-2+t}{e-2} + \binom{e-3+t}{e-2} &\text{ for even $s$.}
    \end{cases}
  \end{equation*}
  By minimality of $Q/I$ one has $t\ge 2$ and hence $\mu(I) \ge
  \binom{e}{e-2} = \frac{e(e-1)}{2} > e$.
\end{rmk}

For artinian Gorenstein local rings of embedding dimension 4---and
more generally for Gorenstein local rings of codepth 4---there is a
classification based on multiplication in Koszul homology. It
predates the classification of local rings of codepth
$3$ and was achieved by Kustin and Miller~\cite{ARKMMl85}; for simplicity we
refer here to Avramov's exposition in \cite{LLA89a}. In addition to the class of
complete intersections, the classification scheme has three classes
one of which is called \textbf{GGO} in \cite{LLA89a}.

\begin{prp}
  \label{prp:ggo}
  An artinian Gorenstein local ring $(R,\m)$ of embedding dimension $4$ and
  socle degree $s$ belongs to the class {\rm \bf GGO} if and only if
  $R/\m^s$ is Golod.
\end{prp}

\begin{prf*}
  Let $k$ denote the residue field of $R$. The Poincar\'e series of
  $R$ has the form $(1+z)^4/d(z)$, where the polynomial $d(z)$ depends
  on the class of $R$ as proved by Jacobsson, Kustin, and
  Miller~\cite{JKM-85}. By \eqref{LA} one
  has
  \begin{equation*}
    \Poin[R/\m^s]{k} = \frac{(1+z)^4}{d(z) - z^2(1+z)^4}\;.    
  \end{equation*}
  For $R/\m^s$ to be Golod, the denominator $D(z) = d(z) - z^2(1+z)^4$
  has to be a polynomial of degree $5$, see \eqref{golod}, but if $R$ is
  not of class \textbf{GGO}, then
  $d(z)$ and hence $D(z)$ has degree at least $7$; see
  \thmcite[(3.5)]{LLA89a}.

  It remains to prove that $R/\m^s$ is Golod if $R$ is of class
  \textbf{GGO}. For a ring $R$ of this class, one gets from
  \thmcite[(3.5)]{LLA89a} the expression
  \begin{align*}
    D(z) &= (1+z)^2(1-2z-(h-3)z^2-2z^3+z^4) -z^2(1+z)^4\\
    &= (1+z)^2(1-2z-(h-2)z^2-4z^3)\\
    &= 1 - (h+1)z^2 -2(h+1)z^3 - (h+6)z^4 -4z^5\:,
  \end{align*}
  where $h$ denotes the minimal number of generators of the defining
  ideal in a minimal Cohen presentation of $R$ (in \cite{LLA89a} this
  number is called $l+1$). Set $h_j = \dim[k]{\H[j]{\Kz^R}}$; as $R$
  is Gorenstein one has $h_0=1=h_4$, $h_1=h=h_3$ and $h_2=2h-2$. By
  \thmcite[1]{GLLLLA78} there is an isomorphism of $k$-algebras
  \begin{equation*}\textstyle
    \H{K^{R/\m^s}} \is \H{K^R}/\H[4]{K^R} \ \ltimes 
    (\mathsf{\Sigma} \bigwedge k^4)/(\mathsf{\Sigma}\bigwedge^4 k^4)\:.
  \end{equation*}
  In particular, one gets
  \begin{align*}
    \dim[k]{\H[1]{\Kz^{R/\m^s}}} &= h_1+1 = h+1\\
    \dim[k]{\H[2]{\Kz^{R/\m^s}}} &= h_2+4 = 2(h+1)\\
    \dim[k]{\H[3]{\Kz^{R/\m^s}}} &= h_3+6 = h+6\\
    \dim[k]{\H[4]{\Kz^{R/\m^s}}} &= 4\:,
  \end{align*}
  and comparison of the expression for $D(z)$ to \eqref{golod}
  shows that $R/\m^s$ is Golod.
\end{prf*}

\begin{rmk}
  \label{rmk:frame}
  Let $(R,\m,k)$ be an artinian Gorenstein local ring of embedding
  dimension $4$ and socle degree $s$. In the terminology of
  \cite{LLA89a}, $R$ is complete intersection or of class
  \textbf{GGO}, \textbf{GT}, or $\mathbf{GH}(p)$; here the parameter
  $p$ is between $1$ and $h-1$ where $h$, as in the proof above, is
  the rank of $\H[1]{K^R}$ or, equivalently, the minimal number of
  generators of the defining ideal in a minimal Cohen presentation of
  $R$.

  Examples of rings of class \textbf{GT} and $\mathbf{GH}$ are
  provided in \prpcite[(2.7) and (2.8)]{ARKMMl85}; they are examples
  of local Gorenstein rings that are not complete intersection and
  have $R/\m^s$ not Golod, compare \prpref{ci}.

  If $R$ has socle degree $3$ and $R$ has an exact zero divisor, then
  \cite{IBHLMS11} yields
  \begin{equation*}
    \Poin[R]{k} = \frac{1}{1-4z+4z^2-z^3} 
    =   \frac{(1+z)^4}{(1+z)^2(1-2z-3z^2+3z^3+2z^4-z^5)}\:,  
  \end{equation*}
  which by \thmcite[3.5]{LLA89a} identifies $R$ as being of class
  $\mathbf{GH}(5)$ where $5=h-1$.

  \exaref{gt} exhibits a concrete local ring $R$ of socle degree $3$ that is
  not complete intersection and does not have an exact zero divisor
  but such that $R/\m^3$ is not Golod, compare \prpref{ezd}.
\end{rmk}

\begin{exa}
  \label{exa:gt}
  Let $k$ be a field and set $Q=\pows{w, x,y,z}$. The ideal
  \begin{equation*}
    I= (w^2+xy, wx+xz, wz, y^2+xz, yz, z^2, x^3+x^2z)
  \end{equation*}
  defines a $k$-algebra with basis
  \begin{equation*}
    1\,;\quad w,x,y,z\,;\quad wy, x^2, xy, xz\,;\quad x^2z\:;
  \end{equation*}
  in particular, it has socle degree $3$. Proceeding as in
  \cite{LWCOVl14a} one can use \textsc{Macaulay~2} \cite{M2} to verify
  that $Q/I$ is a Gorenstein ring of class \textbf{GT}.
\end{exa}

\section*{Acknowledgment}
\noindent
We thank Pedro Macias Marques for discussions related to the material
in this paper, and we thank the anonymous referee
for pertinent questions that prompted us to add \secref{4}.

\bibliographystyle{amsplain-nodash}

\def\soft#1{\leavevmode\setbox0=\hbox{h}\dimen7=\ht0\advance \dimen7
  by-1ex\relax\if t#1\relax\rlap{\raise.6\dimen7
  \hbox{\kern.3ex\char'47}}#1\relax\else\if T#1\relax
  \rlap{\raise.5\dimen7\hbox{\kern1.3ex\char'47}}#1\relax \else\if
  d#1\relax\rlap{\raise.5\dimen7\hbox{\kern.9ex \char'47}}#1\relax\else\if
  D#1\relax\rlap{\raise.5\dimen7 \hbox{\kern1.4ex\char'47}}#1\relax\else\if
  l#1\relax \rlap{\raise.5\dimen7\hbox{\kern.4ex\char'47}}#1\relax \else\if
  L#1\relax\rlap{\raise.5\dimen7\hbox{\kern.7ex
  \char'47}}#1\relax\else\message{accent \string\soft \space #1 not
  defined!}#1\relax\fi\fi\fi\fi\fi\fi}
  \providecommand{\MR}[1]{\mbox{\href{http://www.ams.org/mathscinet-getitem?mr=#1}{#1}}}
  \renewcommand{\MR}[1]{\mbox{\href{http://www.ams.org/mathscinet-getitem?mr=#1}{#1}}}
  \providecommand{\arxiv}[2][AC]{\mbox{\href{http://arxiv.org/abs/#2}{\sf
  arXiv:#2 [math.#1]}}} \def\cprime{$'$}
\providecommand{\bysame}{\leavevmode\hbox to3em{\hrulefill}\thinspace}
\providecommand{\MR}{\relax\ifhmode\unskip\space\fi MR }
\providecommand{\MRhref}[2]{%
  \href{http://www.ams.org/mathscinet-getitem?mr=#1}{#2}
}
\providecommand{\href}[2]{#2}

\end{document}